\documentclass[12pt]{amsart}
\usepackage{amssymb}
\usepackage[mathscr]{eucal}
\usepackage{amscd}
\usepackage{amsfonts}
\usepackage{amsmath, amsthm, amssymb}
\usepackage{latexsym}
\usepackage{hyperref}
\usepackage{float}
\usepackage{graphicx}
 \theoremstyle{plain}

\theoremstyle{definition}

\theoremstyle{remark}

\allowdisplaybreaks

\begin{document}
\title[Solution of a Nonlinear Integral Equation]
{Solution of a Nonlinear Integral Equation via New Fixed Point
Iteration Process}

\author[Garodia and Uddin]{ Chanchal Garodia and Izhar Uddin}
\thanks{$^{\ast}$Corresponding author: Izhar Uddin, Department of mathematics, Jamia Millia Islamia, New Delhi-110025, India. Email address:
 {\rm izharuddin1@jmi.ac.in}}
\maketitle
\begin{center}
{\footnotesize Department of Mathematics, Jamia Millia Islamia, New Delhi-110025, India.\\

Email addresses: c.garodia85@gmail.com and izharuddin1@jmi.ac.in }
\end{center}

{\footnotesize \noindent {\bf Abstract.} In this paper, we introduce
a new three-step iteration process in Banach space and prove
convergence results for approximating fixed points for nonexpansive
mappings. Also, we show that the newly introduced iteration process
converges faster than a number of existing iteration processes.
Further, we discuss about the solution of mixed type
Volterra-Fredholm functional nonlinear integral equation.}

 \vskip0.5cm \noindent {\bf Keywords}: Fixed point, Strong convergence, Weak convergence, Nonexpansive mapping, Banach space.

\noindent {\bf AMS Subject Classification}: 47H09, 47H10.

\section{INTRODUCTION}
The study of fixed point theorems and their applications though
initiated long ago, still continues to be a highly challenging yet
useful area of investigation in Topology and Nonlinear functional
analysis. The existing literature on fixed point theory is
extensive. The study of fixed points lies within several domain
which include: Topology, Algebraic topology, Nonlinear operators,
Differential equations (Ordinary and Partial both) and Functional
analysis. Unlike many domains of pure mathematics, the fixed point
theory has numerous applications in various fields within as well as
beyond the mathematics namely: Approximation theory, Successive
approximation, Integral equations, Game theory, Optimal control,
Optimization, Economics and several others.\\
A wide range of problems of applied sciences and engineering are
usually formulated as functional equations. Such equations can be
written in the form of fixed point equations. Operator equations
representing phenomena occurring in different fields, such as steady
state temperature distribution, chemical reactions, neutron
transport theory, economic theories and epidemics, often require
appropriate and adequate solutions. Thus, the aim of finding
solution to these equations is to locate the fixed point and
approximate its value. However, once we ensure the existence of a
fixed point of some mapping, then it is always desirable to develop
such methods which can be efficiently used to approximate that fixed
point. Iterative process are one of the fundamental tool that can be
used to locate a fixed point. In the last few years, various authors
have introduced numerous iterative schemes which have been utilized
widely to approximate the fixed points of operators. Banach
contraction theorem \cite{SB} which is one of the most widely and
extensively utilized result use Picard iteration process for
locating the fixed point. Owing to the importance of iteration
processes, many new iteration schemes have been obtained in the last
few years and the prime focus of researchers is to obtain the
iteration schemes which converges at a faster rate than the existing
schemes. Some of the well known iteration processes are Mann
iteration \cite{WM}, Ishikawa iteration \cite{SI}, Halpern iteration
\cite{BHN}, Noor iteration \cite{MAN}, Agarwal et al. iteration
\cite{RPA}, SP iteration \cite{SP},
Normal-S iteration \cite{NS} and Abbas and Nazir iteration \cite{AT}.\\
 Let $J$ be a nonempty closed convex subset of a uniformly convex Banach space $G$. A mapping $S : J\rightarrow J$
is said to be nonexpansive if $\|Su - Sv\| \leq \|u - v\|$ for all
$u, v \in J.$ A point $q \in J$ is said to be a fixed point of $S$
if $Sq = q$. We will use $F(S)$ to denote the set of fixed points of $S$.\\
Recently, Thakur et al. \cite{TTP} introduced the following
iteration, where a sequence $\{w_n\}$ is constructed from arbitrary
$w_1 \in J$ by:
$$ {\begin{cases}

u_n = (1-\delta_n)w_n + \delta_nSw_n \cr
 v_n = S((1-\zeta_n)w_n + \zeta_nu_n) \cr
w_{n+1} = Sv_n
\end{cases}} \eqno {(1.1)}$$
 where $\{\delta_n\}$ and $\{\zeta_n\}$ are in $(0, 1).$ They proved that their process converges faster than Abbas and Nazir iteration \cite{AT}.\\
Motivated and inspired by the research going on in this direction,
we introduce a new iteration process for approximating fixed point
of a nonexpansive mapping, where the sequence $\{c_n\}$ is generated
iteratively by $c_1 \in J$ and
$$ {\begin{cases}

a_n =Sc_n \cr
 b_n = (1-\delta_n)a_n + \delta_nSa_n \cr
c_{n+1} = Sb_n
\end{cases}} \eqno {(1.2)}$$
for each $n \in \mathbb{N}$ and $\{\delta_n\}$ is a sequence in $(0,
1).$\\
The purpose of this paper is to prove the convergence of newly
defined iteration process $(1.2)$ for nonexpansive mappings. We
prove that iteration process $(1.2)$ converges faster than Thakur et
al. iteration $(1.1)$ which is faster than Picard, Mann, Ishikawa,
Noor, Agarwal et al. and Abbas and Nazir iteration processes. We
also present numerical example to compare the convergence of $(1.2)$
with Agarwal et al., Abbas and Nazir, Noor and Thakur et al.
iteration schemes. In addition to this, we show that newly
introduced iteration scheme $(1.2)$ converges strongly to a solution
of nonlinear integral equation.

\section{PRELIMINARIES}
We begin by recalling some known Results and Definitions which will
be frequently used through out the text.\\
{\noindent\bf{Definition 2.1.}} A Banach space $G$ is said to be
uniformly convex if for each $\alpha \in (0, 2]$ there is a $\beta
> 0$ such that for $a, b \in G$ with $\|a\| \leq 1$, $\|b\| \leq 1$
and $\|a - b\| > \alpha$, we have $$\left\|\frac{a + b}{2}\right\| <
\beta.$$
 {\noindent\bf{Definition 2.2.}} A Banach space $G$ is
said to satisfy the Opial's condition if for any sequence $\{a_n\}$
in $G$ which converges weakly to $a \in G$ i.e. $a_n \rightharpoonup
a$ implies that $$\limsup_{n\to\infty} \|a_n - a\| <
\limsup_{n\to\infty} \|a_n - b\|$$ for all $b \in G$ with $b
\neq a.$\\
A mapping $S : J \rightarrow G$ is demiclosed at $a \in G$ if for
each sequence $\{a_n\}$ in $J$ and each $b \in G$, $a_n
\rightharpoonup b$ and $Sa_n \rightarrow a$ imply that $b \in J$ and
$Sb = a.$\\
The following definitions about the rate of convergenve were given
by Berinde \cite{VB}.\\
{\noindent\bf{Definition 2.3.}} Let $\{r_n\}$ and $\{p_n\}$ be two
real sequences converging to $r$ and $p$ respectively. Then,
$\{r_n\}$ converges faster then $\{p_n\}$ if
$\lim\limits_{n\to\infty}\frac{\|r_n - r\|}{\|p_n - p\|} = 0.$\\
{\noindent\bf{Definition 2.4.}} Let $\{y_n\}$ and $\{z_n\}$ be two
fixed point iteration processes converging to the same fixed point
$q$. If $\{r_n\}$ and $\{p_n\}$ are two sequences of positive
numbers converging to zero such that $\|y_n - q\| \leq r_n$ and
$\|z_n - q\| \leq p_n$ for all $n \geq 1$, then we say that
$\{y_n\}$ converges faster than $\{z_n\}$ to $q$ if $\{r_n\}$
converges faster then $\{p_n\}$.\\
Next, we list two Lemmas which will be useful in our subsequent
discussion.\\
{\noindent\bf{Lemma 2.1.}} (\cite{KKT}) Let $J$ be a nonempty closed
convex subset of a uniformly convex Banach space $G$ and $S$ a
nonexpansive mapping on $J.$ Then, $I - S$ is demiclosed at zero.\\
{\noindent\bf{Lemma 2.2.}} (\cite{JS}) Let $G$ be a uniformly convex
Banach space and $\{u_n\}$ be any sequence such that $0 < w \leq u_n
\leq v < 1$ for some $w, v \in \mathbb{R}$ and for all $n \geq 1.$
Let $\{p_n\}$ and $\{r_n\}$ be any two sequences of $G$ such that
$\limsup\limits_{n\to\infty}\|p_n\| \leq q$,
$\limsup\limits_{n\to\infty}\|r_n\| \leq q$ and
$\limsup\limits_{n\to\infty}\|u_np_n + (1-u_n)r_n\| = q$ for some $q
\geq 0.$ Then, $\lim\limits_{n\to\infty}\|p_n - r_n\| = 0.$\\

\section{RESULTS}
In this section, first we show that our iteration scheme $(1.2)$
converges faster than the iteration of Thakur et al. $(1.1).$\\
%
%
%
%
{\noindent\bf{Theorem 3.1.}} Let $J$ be a nonempty closed convex
subset of a Banach space $G$ and $S: J \rightarrow J$ be a
contraction mapping with contraction factor $\xi \in (0, 1)$ such
that $F(S) \neq \Phi.$ 
If $\{c_n\}$ is a sequence defined by $(1.2)$, then $\{c_n\}$
converges faster than the iteration scheme of Thakur et al. $(1.1).$\\
{\noindent\bf{Proof.}} From $(1.2)$, for any $q \in F(S)$,
$$\| a_n - q\| = \| Sc_n - q\| \leq  \xi \| c_n - q\|$$
and
$$\begin{array}{lll}
\| b_n - q \| &=& \| (1-\delta_n)a_n + \delta_nSa_n - q\| \\
&\leq& (1-\delta_n)\|a_n - q\| + \xi\delta_n\|a_n - q\|\\
&=& (1 - (1 - \xi)\delta_n)\|a_n - q\|\\
&\leq& \xi(1 - (1 - \xi)\delta_n)\|c_n - q\|.
\end{array}$$
So,
$$\begin{array}{lll}
\|c_{n+1} - q\| &=& \|Sb_n - q\| \\
&\leq& \xi\|b_n - q\|\\
&\leq& \xi^2 (1 - (1 - \xi)\delta_n)\|c_n - q\|\\
 &\leq& \xi^{2n} (1 - (1 - \xi)\delta)^n\|c_1 - q\|.
\end{array}$$
Now, using $(1.1)$, we obtain
$$\begin{array}{lll}
\|u_n - q\| &=& \|(1-\delta_n)w_n + \delta_nSw_n - q\| \\
&\leq&(1-\delta_n)\|w_n - q\| + \delta_n\|Sw_n - q\|\\
&\leq&(1-\delta_n)\|w_n - q\| + \xi\delta_n\|w_n - q\|\\
&=&(1-(1-\xi)\delta_n)\|w_n - q\|
\end{array}$$
and
$$\begin{array}{lll}
\|v_n - q\| &=& \|S((1-\zeta_n)w_n + \zeta_n u_n)- q\| \\
&\leq& \xi\|(1-\zeta_n)w_n + \zeta_n u_n- q\| \\
&\leq&\xi((1-\zeta_n)\|w_n - q\| + \zeta_n\|u_n - q\|)\\
&\leq&\xi((1-\zeta_n)\|w_n - q\| + \zeta_n(1-(1-\xi)\delta_n)\|w_n -
q\|)\\
&=&\xi(1-(1 - \xi)\delta_n\zeta_n)\|w_n - q\|.

\end{array}$$
Thus,
$$\begin{array}{lll}
\|w_{n+1} - q\| &=& \|Sv_n - q\|\\
&\leq& \xi\|v_n - q\|\\
&\leq& \xi^2(1-(1 - \xi)\delta_n\zeta_n)\|w_n - q\|\\
&\leq& \xi^{2n}(1-(1 - \xi)\delta\zeta)^n\|w_1 - q\|

\end{array}$$
Let $$r_n = \xi^{2n}(1-(1 - \xi)\delta\zeta)^n\|w_1 - q\|$$ and
$$p_n = \xi^{2n} (1 - (1 - \xi)\alpha)^n\|c_1 - q\|.$$ Then,
$$\begin{array}{lll}  \frac{p_n}{r_n} &=& \frac{\xi^{2n} (1 - (1 -
\xi)\delta)^n\|c_1 - q\|}{\xi^{2n}(1-(1 - \xi)\delta\zeta)^n\|w_1 - q\|}\\
&=& \frac{(1 - (1 -
\xi)\delta)^n\|c_1 - q\|}{(1-(1 - \xi)\delta\zeta)^n\|w_1 - q\|}\\
&\rightarrow& ~~0 ~~~~~~~~~~~as ~~~~~~~~~~~~ n \rightarrow \infty.
\end{array}$$
Thus, $\{c_n\}$ converges faster than $\{w_n\}.$\\
{\noindent\bf{Lemma 3.1.}} Let $J$ be a nonempty closed convex
subset of a Banach space $G$ and $S: J \rightarrow J$ be a
nonexpansive mapping with $F(S) \neq \Phi.$ Let $\{c_n\}$ be the
iterative sequence defined by the iteration process $(1.2)$. Then,\\
(i) $\lim\limits_{n\to\infty} \|c_n - q\|$ exists for all $q \in F(S),$\\
(ii) $\lim\limits_{n\to\infty} \|Sc_n - c_n\| = 0.$\\
 {\noindent\bf{Proof.}} (i) Let $q \in F(S)$. Then, using $(1.2)$ we obtain
  $$\| a_n - q\| = \| Sc_n - q\| \leq \| c_n - q\| \eqno {(3.1)}$$
  and
$$\begin{array}{lll}
\| b_n - q \| &=& \| (1-\delta_n)a_n + \delta_nSa_n - q\| \\
&\leq& (1-\delta_n)\|a_n - q\| + \delta_n\|Sa_n - q\|\\
&\leq& (1-\delta_n)\|a_n - q\| + \delta_n\|a_n - q\| \\
&=&\|a_n - q\|\\
 &\leq& \|c_n - q\|.
\end{array} \eqno {(3.2)}$$
Using $(3.1)$ and $(3.2)$, we get
$$\begin{array}{lll}
\|c_{n+1} - q\| &=& \|Sb_n - q\| \\
&\leq& \|b_n - q\|\\
&\leq&\|c_n - q\|.
\end{array}$$
Thus, $\lim\limits_{n\to\infty} \|c_n - q\|$ exists for all $q \in
F(S).$\\
(ii) Let $\lim\limits_{n\to\infty} \|c_n - q\| = \kappa.$ \\
From $(3.1)$ and $(3.2)$, we have
$$\limsup\limits_{n\to\infty} \|b_n - q\| \leq \kappa  \eqno {(3.3)} $$ and
$$\limsup\limits_{n\to\infty} \|a_n - q\| \leq \kappa.  \eqno {(3.4)}$$
Now, $$\kappa = \lim\limits_{n\to\infty} \|c_{n+1} - q\| =
\lim\limits_{n\to\infty} \|Sb_n - q\|,$$ and $$\|Sb_n - q\| \leq
\|b_n - q\|.$$ So, $$\kappa \leq \liminf\limits_{n\to\infty} \|b_n -
q\|$$ which along with $(3.3)$ implies
$$\lim\limits_{n\to\infty} \|b_n - q\| = \kappa.  \eqno {(3.5)}$$
As, $$\|Sa_n - q\| \leq \|a_n - q\|$$
from $(3.4)$, we obtain
$$\limsup\limits_{n\to\infty} \|Sa_n - q\| \leq \kappa .  \eqno {(3.6)}$$
Consider,
$$\begin{array}{lll}
\lim\limits_{n\to\infty} \|b_n - q\| &=&
\lim\limits_{n\to\infty}\|(1-\delta_n)a_n + \delta_nSa_n - q\|\\
&=& \lim\limits_{n\to\infty}\|(1-\delta_n)(a_n - q) + \delta_n(Sa_n
- q)\|.
\end{array}$$
Using Lemma 2.2, from $(3.4)$, $(3.5)$ and $(3.6)$, we get
$$\lim\limits_{n\to\infty}\|a_n - Sa_n\| = 0.  \eqno {(3.7)}$$
Now, consider
$$\begin{array}{lll}
\|b_n - Sa_n\| &=& \|(1-\delta_n)a_n + \delta_nSa_n - Sa_n\|\\
&=&\|(1-\delta_n)(a_n - Sa_n)\|,
\end{array}$$ which on using $(3.7)$ gives
$$\lim\limits_{n\to\infty}\|b_n - Sa_n\| = 0.  \eqno {(3.8)}$$
Since, $$\|a_n - b_n\| \leq \|a_n - Sa_n\| + \|Sa_n - b_n\|,$$ this
together with $(3.7)$ and $(3.8)$ yields that
$$\lim\limits_{n\to\infty}\|a_n - b_n\| = 0.  \eqno {(3.9)}$$
Now, using $(3.8)$ and $(3.9)$, we have
$$\begin{array}{lll}
\|Sc_{n+1} - c_{n+1}\| &=& \|Sc_{n+1} - Sb_n\|\\
&\leq& \|c_{n+1} - b_n\|\\
&=& \|Sb_n - b_n\|\\
&=& \|Sb_n - Sa_n + Sa_n - b_n\|\\
&\leq& \|Sb_n - Sa_n\| + \|Sa_n - b_n\|\\
&\leq& \|b_n - a_n\| + \|Sa_n - b_n\|
\end{array}$$
Hence, $$\lim\limits_{n\to\infty}\|Sc_{n} - c_{n}\| = 0.$$ Now, we
prove the weak convergence of iteration process $(1.2).$\\
{\noindent\bf{Theorem 3.2.}} Let $J$ be a nonempty closed convex
subset of a uniformly convex Banach space $G$ which satisfies the
Opial's condition and $S: J \rightarrow J$ be a nonexpansive mapping
with $F(S) \neq \Phi.$ If $\{c_n\}$ is the iterative sequence
defined by the iteration process $(1.2)$, then $\{c_n\}$ converges
weakly to a fixed point of $S$.\\
{\noindent\bf{Proof.}} Let $q\in F(S)$. Then, from Lemma 3.1
$\lim\limits_{n\to\infty}\|c_n - q\|$ exists. In order to show the
weak convergence of the iteration process $(1.2)$ to a fixed point
of $S$, we will prove that $\{c_n\}$ has a unique weak subsequential
limit in $F(S).$ For this, let $\{c_{n_\upsilon}\}$ and
$\{c_{n_\rho}\}$ be two subsequences of $\{c_n\}$ which converges
weakly to $u$ and $v$ respectively. By Lemma 3.1, we have
$\lim\limits_{n\to\infty}\|Sc_{n} - c_{n}\| = 0$ and using the Lemma
2.1, we have $I - S$ is demiclosed at zero. So $u, v \in F(S)$.\\
Next, we show the uniqueness. Since $u, v \in F(S)$, so
$\lim\limits_{n\to\infty}\|c_n - u\|$ and
$\lim\limits_{n\to\infty}\|c_n - v\|$ exists. Let $u \neq v$. Then,
by Opial's condition, we obtain
$$\begin{array}{lll}
\lim\limits_{n\to\infty}\|c_n - u\| &=&
\lim\limits_{n\to\infty}\|c_{n_\upsilon} - u\|\\ &<&
\lim\limits_{n\to\infty}\|c_{n_\upsilon} - v\| \\&=&
\lim\limits_{n\to\infty}\|c_n - v\|\\ &=&
\lim\limits_{n\to\infty}\|c_{n_\rho} - v\|\\ &<&
\lim\limits_{n\to\infty}\|c_{n_\rho} - u\| \\&=&
\lim\limits_{n\to\infty}\|c_n - u\|
\end{array}$$
which is a contradiction, so $u = v.$ Thus, $\{c_n\}$ converges
weakly to a fixed point of $S$.\\
Next, we establish some strong convergence results for iteration process $(1.2).$\\
{\noindent\bf{Theorem 3.3.}} Let $J$ be a nonempty closed convex
subset of a uniformly convex Banach space $G$ and $S: J \rightarrow
J$ be a nonexpansive mapping with $F(S) \neq \Phi.$ If $\{c_n\}$ is
the iterative sequence defined by the iteration process $(1.2)$,
then $\{c_n\}$ converges to a point of $F(S)$ if and only if $\liminf\limits_{n\to\infty}d(c_n, F(S)) = 0.$\\
{\noindent\bf{Proof.}} If the sequence $\{c_n\}$ converges to a
point $q \in F(S)$, then it is obvious that $\liminf\limits_{n\to\infty} d(c_n, F(S))=0.$\\
For converse part, assume that $\liminf\limits_{n\to\infty} d(c_n,
F(S))=0.$ From Lemma 3.1, we have $\lim\limits_{n\to\infty}\|c_n -
q\|$ exists for all $q \in F(S)$, which gives
$$\|c_{n+1} - q\| \leq \|c_n - q\|~\mbox{for~ any} ~q\in F(S)$$
which yields
$$d(c_{n+1}, F(S)) \leq d(c_{n}, F(s)).$$
Thus, $\{d(c_n, F(S))\}$ forms a decreasing sequence which is
bounded below by zero as well, so we get that
$\lim\limits_{n\to\infty} d(c_n, F(S))$
 exists. As, $\liminf\limits_{n\to\infty} d(c_n, F(S))=0$ so $\lim\limits_{n\to\infty} d(c_n, F(S))=0.$\\
 Now, we prove that $\{c_n\}$ is a cauchy sequence in $J$. Let $\epsilon > 0$ be arbitrarily chosen.
  Since $\liminf\limits_{n\to\infty} d(c_n, F(S))=0$, there
exists $n_0$ such that for all $n \geq n_0$, we have
$$d(c_n, F(S) ) < \frac{\epsilon}{4}.$$
In particular,
$$\inf \{ \|c_{n_0} - q\|: q\in F(S) \}
< \frac{\epsilon}{4},$$ so there must exist a $\theta \in F(S)$ such
that
$$ \|c_{n_0} - \theta\|< \frac{\epsilon}{2}.$$
  Thus, for $m, n \geq n_0$, we have
$$\| c_{n+m} - c_n\| \leq  \| c_{n+m} - \theta\| + \|c_n - \theta\| <  2\|c_{n_0} - \theta\| <
2\frac{\epsilon}{2} = {\epsilon}$$ which shows that $\{c_n\}$ is a
cauchy sequence. Since $J$ is a closed subset of a Banach space $G$,
therefore $\{c_n\}$ must converge in $J$. Let,
$\lim\limits_{n\to\infty}c_n = r$ for some $r \in J$.

Now, using $\lim\limits_{n\to\infty} \|Sc_n - c_n\|=0$, we get
$$\begin{array}{lll}
 \|r - Sr\| &\leq& \|r - c_n\| + \|c_n - Sc_n\| + \|Sc_n - Sr\|\\
 &\leq& \|r - c_n\| + \|c_n - Sc_n\| + \|c_n - r\|  \end{array}$$
 $$\rightarrow 0~\mbox{as}~ n \rightarrow \infty$$
and hence $r = Sr.$ Thus, $r \in F(S).$ This proves our result.
%


 In (\cite{HFS}), Senter and Dotson gave the condition (A) which states
that a mapping $S: J \rightarrow J$ is said to satisfy the
condition(A) if there exists a nondecreasing function $g:[0,\infty)
\rightarrow [0,\infty)$ with $g(0)=0$ and $g(r)>0$ for all $r\in
(0,\infty)$ such that $\|u - Su\| \geq g(d(u, F(S)))$ for
all $u \in J,$ where $d(u, F(S)) = inf\{\|u - q\| : q \in F(S)\}.$\\

{\noindent\bf{Theorem 3.4.}} Let $J$ be a nonempty closed convex
subset of a uniformly convex Banach space $G$. Let $S: J \rightarrow
J$ be a nonexpansive mapping such that $F(S)\neq \phi$ and $\{c_n\}$
be the sequence defined by $(1.2)$. If $S$ satisfies condition (A),
then $\{c_n\}$ converges strongly to a fixed point of $S$. \\
{\noindent\bf{Proof.}} By Lemma 3.1, $\lim\limits_{n\to\infty} \|c_n
- q\|$ exists and $\|c_{n+1} - q\| \leq \|c_n - q\|$ for all $q \in F(S).$\\
We get $$ \inf\limits_{q \in F(S)}\|c_{n+1} - q\| \leq
\inf\limits_{q \in F(S)} \|c_n - q\|,$$ which yields
$$d(c_{n+1}, F(S))\leq d(c_n, F(S)).$$
This shows that the sequence $\{d(c_n, F(S))\}$ is decreasing and
bounded below, so $\lim\limits_{n\to\infty} d(c_n, F(S))$ exists.\\
Also, by Lemma 3.1 we have $\lim\limits_{n\to\infty} \|c_n - Sc_n\|=0.$\\
It follows from condition (A) that $$\lim\limits_{n\to\infty}
g(d(c_n, F(S)))\leq \lim\limits_{n\to\infty} \|c_n - Sc_n\|=0,$$
so that $\lim\limits_{n\to\infty} g(d(c_n, F(S)))=0.$\\
Since $g$ is a non decreasing function satisfying $g(0) = 0$ and
$g(r)> 0$ for all $r \in (0, \infty)$, therefore
$\lim\limits_{n\to\infty} d(c_n,
F(S))=0.$\\
By Theorem 3.3., the sequence $\{c_n\}$ converges strongly to a
point of $F(S).$
\section{NUMERICAL EXAMPLE}
In this section, we present an example which shows that our
iteration process $(1.2)$ converges faster than Noor, Agarwal et
al., Abbas and Nazir and
Thakur et al. iteration processes.\\
{\noindent\bf{Example}} Let $G = \mathbb{R}$ and $J = [1, 50]$. Let
$S: J \rightarrow J$ be a mapping defined as $S(c) = \sqrt{c^2 - 6c
+ 30}$ for all $c \in J$. Clearly $c = 5$ is the fixed point of $S$.
Set $\delta_n = 0.95$, $\zeta_n = 0.30$ and $\gamma_n = 0.90$ for
all $n \in \mathbb{N}.$ Choose initial value as $40.$ Then, we get
the following table of iteration values:

\begin{table}[H]
\begin{tabular}{ |c|c|c|c|c|c| }
 \hline
Step& Agarwal & Noor& Abbas& Thakur&New iter \\
\hline
1& 40& 40& 40& 40& 40\\
2& 36.514581536& 36.1407358454 &34.6327094201 &33.826119187&
32.0516661514\\
3& 33.0679397292& 32.3319800634& 29.367596341& 27.7952382214&
24.3704817561\\
 4& 29.6685026434& 28.5864001018& 24.2441205249&
21.9755503569 &17.1512394673\\
 5 &26.3275706425& 24.9216680332&
19.3272253806& 16.4912748872& 10.8507668765\\
6& 23.0606740255& 21.3631923337 &14.7312047769& 11.58909835&
6.56090714498\\
 7 &19.8897435795 &17.9486804932 &10.669341116&
7.76797689976& 5.170067389 \\
8 &16.8466375807& 14.7357042537 &7.52753037638 &5.69419487242&
5.01234906914\\
9& 13.9788293647& 11.8131560009 &5.74872981287 &5.11047433488
&5.00085298124\\
 10& 11.3578863208 &9.31239170384& 5.15293698217&
5.01509864089& 5.00005870115\\
 11& 9.08838738905& 7.39232766846&
5.02680130355& 5.002010937& 5.00000403871\\
 12 &7.30211704972&6.14706302858& 5.00453782329& 5.00026687693& 5.00000027786 \\
 13&6.10042750793 &5.48679502226& 5.00076360356& 5.00003540113&
5.00000001912\\
 14& 5.44793732586 &5.19164605164& 5.00012836141&
5.00000469565 &5.00000000132 \\
15& 5.16283265193& 5.07277137427& 5.0000215737 &5.00000062283&
5.00000000009\\
 16& 5.05591482526 &5.02722036618 &5.00000362578&
5.00000008261& 5.00000000001 \\
17& 5.0187690152& 5.01012264162 &5.00000060936& 5.00000001096  &
5.00000000000\\
 18& 5.00624950826 &5.00375610854 &5.00000010241&
5.00000000145& 5.00000000000\\
 19& 5.00207518849& 5.00139259846&
5.00000001721 &5.00000000019& 5.00000000000\\
20& 5.00068844699& 5.00051615633& 5.00000000289 &5.00000000003&
5.00000000000\\
 21& 5.00022832364& 5.00019128789& 5.00000000049  & 5.00000000000&
 5.00000000000\\
22& 5.00007571593& 5.00007088845& 5.00000000008 &5.00000000000&
5.00000000000\\
 23& 5.00002510782& 5.0000262698 &5.00000000001& 5.00000000000& 5.00000000000\\
 24&5.0000083258& 5.00000973499&   5.00000000000& 5.00000000000&  5.00000000000\\
 25& 5.00000276084&5.00000360756&   5.00000000000&  5.00000000000& 5.00000000000\\
 26 & 5.0000009155& 5.00000133688&   5.00000000000& 5.00000000000&
 5.00000000000\\
27&  5.00000030358& 5.00000049541&   5.00000000000&  5.00000000000&
5.00000000000\\
28& 5.00000010067 &5.00000018359&   5.00000000000&  5.00000000000&
5.00000000000\\
 29& 5.00000003338& 5.00000006803&   5.00000000000&  5.00000000000& 5.00000000000\\
 30&5.00000001107& 5.00000002521  & 5.00000000000& 5.00000000000& 5.00000000000\\
\hline
\end{tabular}
\caption{}
\end{table}

 Also, the following graph shows that our iteration process
$(1.2)$ converges faster to $c = 5$ which is a fixed point of $S$.
\begin{figure}[H]
 \frame{\includegraphics[scale=0.52]{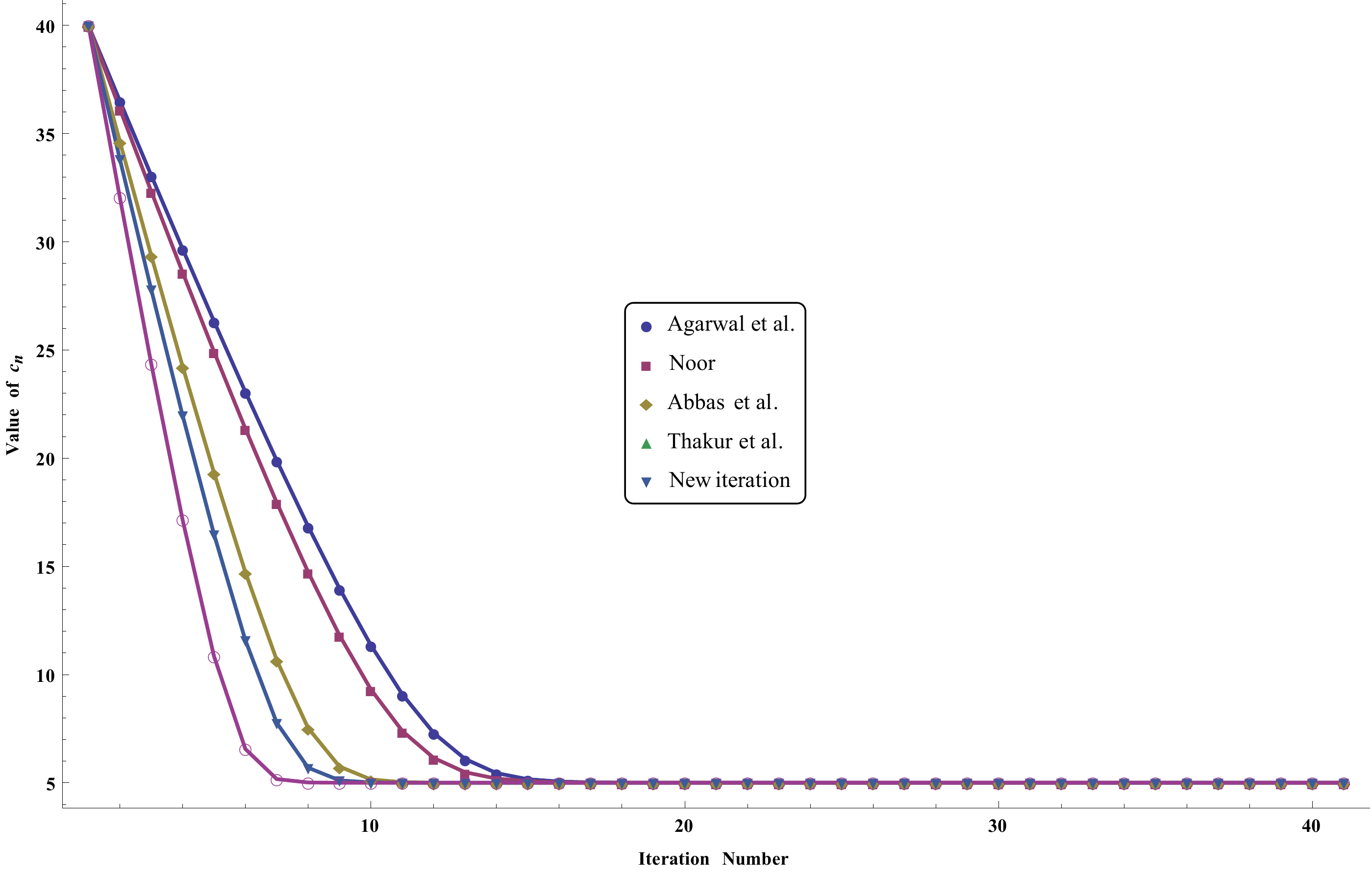}}
\caption{Convergence behaviour of iteration (1.2)}
\end{figure}
Thus, it is evident from the Table 1 as well as Figure 1  that the
newly introduced iteration scheme converges at a much faster rate
than a number of the existing iteration processes.
\section{APPLICATION}
In this section, we will give some application of newly introduced
iteration scheme. We will show that iterative algorithm $(1.2)$
converges strongly to the solution of the following mixed type
Volterra-Fredholm functional nonlinear integral equation which was
discussed in \cite{vol} :
$$x(t) = F \left (t, x(t), \int_{g_1}^{j_1} ... \int_{g_m}^{j_m} K(t, s,
x(s)) ds, \int_{g_1}^{h_1} ... \int_{g_m}^{h_m} H(t, s, x(s)) ds
\right ), \eqno{(5.1)}$$ where $[g_1; h_1] \times ... \times [g_m;
h_m]$ is an interval in $\mathbb{R}^{m},$ K, H : $[g_1; h_1] \times
... \times [g_m; h_m] \times [g_1; h_1] \times ... \times [g_m; h_m]
\times \mathbb{R} \rightarrow \mathbb{R}$ continuous functions and
$F: [g_1; h_1] \times ... \times [g_m; h_m] \times \mathbb{R}^{3}
\rightarrow
\mathbb{R}.$\\
First, we recall the following theorem which will be fruitful in our
subsequent discussion.\\
{\noindent\bf{Theorem 5.1.}} (\cite{vol}) We suppose that the following conditions are satisfied:\\
 $(A_1)$ $K, H \in C([g_1; h_1] \times ... \times [g_m; h_m]
\times [g_1; h_1] \times ... \times [g_m; h_m] \times
\mathbb{R})$;\\
$(A_2)$ $F \in C([g_1; h_1] \times ... \times [g_m; h_m] \times
\mathbb{R}^{3});$\\
$(A_3)$ there exist nonnegative constants $\alpha, \beta$ and
$\gamma$ such that
$$|F(t, u_1, v_1, w_1) - F(t, u_2, v_2, w_2)| \leq \alpha |u_1 - u_2|
+ \beta |v_1 - v_2| + \gamma |w_1 - w_2|,$$ for all $ t \in [g_1;
h_1] \times ... \times [g_m; h_m],$ $u_i, v_i,
w_i \in \mathbb{R}, i = 1,2;$\\
$(A_4)$ there exist nonnegative constants $L_K$ and $L_H$ such that
$$|K(t, s, u) - K(t, s, v)| \leq L_K |u - v|,$$
$$|H(t, s, u) - H(t, s, v)| \leq L_H |u - v|,$$ for all $t, s \in
[g_1; h_1] \times ... \times [g_m; h_m],$ $u, v
\in \mathbb{R};$\\
$(A_5)$ $\alpha +( \beta L_K + \gamma L_H) (h_1 - g_1) ...(h_m -
g_m) < 1.$\\
Then $(5.1)$ has a unique solution $p \in C([g_1; h_1] \times
... \times [g_m; h_m]).$\\
Now, we prove our main theorem.\\
{\noindent\bf{Theorem 5.2.}} One opines that all conditions $(A_1) -
(A_5)$ in Theorem 5.1 are performed. Let $\{\delta_n\} \subset [0,
1]$ be a real sequence satisfying
$\sum\limits_{n=1}^{\infty}\delta_n = \infty.$ Then, $(5.1)$ has a
unique solution say $p$ in $C([g_1; h_1] \times ... \times [g_m;
h_m])$ and iteration scheme $(1.2)$ converges to $p.$\\
{\noindent\bf{Proof.}} We consider the Banach space $ G = C([g_1;
h_1] \times ... \times [g_m; h_m], {\|.\|}_{C}),$ where
${\|.\|}_{C}$ is Chebyshev's norm. Let $\{c_n\}$ be an iterative
sequence generated by our iteration scheme for the operator $A: G
\rightarrow G$ defined by
$$A(c)(t) =  F \Big(t, c(t), \int_{g_1}^{j_1} ... \int_{g_m}^{j_m}
K(t, s, c(s)) ds, \int_{g_1}^{h_1} ... \int_{g_m}^{h_m} H(t, s,
c(s)) ds\Big ). \eqno{(5.2)}$$ We will show that $c_n \rightarrow p$
as $n \rightarrow \infty.$\\
From $(1.2)$, $(5.1)$, $(5.2)$ and assumptions $(A_1) - (A_4),$ we
get that
\begin{equation*}
\begin{aligned}
\| c_{n+1} - p\| & = \| Ab_n - p\|  \\
                 & = | A(b_n)(t) - A(p)(t)| \\
                  & = \Big | F \Big (t, b_n(t), \int_{g_1}^{j_1} ... \int_{g_m}^{j_m} K(t,
s, b_n(s)) ds, \int_{g_1}^{h_1} ... \int_{g_m}^{h_m} H(t, s, b_n(s))ds \Big ) \\
& \hspace{0.5cm} -  F \Big (t, p(t), \int_{g_1}^{j_1} ...
\int_{g_m}^{j_m} K(t, s, p(s)) ds, \int_{g_1}^{h_1} ...
\int_{g_m}^{h_m} H(t, s, p(s)) ds \Big ) \Big | \\
& \leq  \alpha |b_n(t) - p(t)| +
  \beta \Big |\int_{g_1}^{j_1} ... \int_{g_m}^{j_m} K(t, s,
b_n(s))ds - \int_{g_1}^{j_1} ... \int_{g_m}^{j_m} K(t, s, p(s))
ds \Big|\\
& \hspace{0.5cm} +
 \gamma \Big |\int_{g_1}^{h_1} ... \int_{g_m}^{h_m} H(t, s, b_n(s)) ds
- \int_{g_1}^{h_1} ... \int_{g_m}^{h_m} H(t, s, p(s)) ds \Big| \\
& \leq \alpha |b_n(t) - p(t)| + \beta \int_{g_1}^{j_1} ...
\int_{g_m}^{j_m} |K(t, s, b_n(s)) - K(t, s, p(s))|ds\\
& \hspace{0.5cm}  + \gamma \int_{g_1}^{h_1} ... \int_{g_m}^{h_m}
|H(t, s, b_n(s))-H(t, s, p(s))|ds\\
& \leq \alpha |b_n(t) - p(t)| + \beta \int_{g_1}^{j_1} ... \int_{g_m}^{j_m} L_K |b_n(s) - p(s)| ds \\
& \hspace{0.5cm} + \gamma \int_{g_1}^{h_1} ... \int_{g_m}^{h_m} L_H |b_n(s) - p(s)|ds\\
& \leq \Big [ \alpha +(\beta L_K + \gamma L_H)
\prod\limits_{i=1}^{m}(h_i- g_i) \Big] \|b_n - p\|.
 \hspace{5cm}(5.3)
\end{aligned}
\end{equation*}

\begin{equation*}
\begin{aligned}
\|b_n - p\| &= \| (1-\delta_n)a_n + \delta_nAa_n - p\|\\
& \leq (1 - \delta_n)| a_n(t) - p(t)| + \delta_n |A(a_n)(t)- A(p)(t)|\\
&= (1 - \delta_n)| a_n(t) - p(t)|\\
& \hspace{0.5cm} + \delta_n \Big | F \Big(t, a_n(t),
\int_{g_1}^{j_1} ... \int_{g_m}^{j_m} K(t, s, a_n(s)) ds, \int_{g_1}^{h_1} ... \int_{g_m}^{h_m} H(t, s, a_n(s)) ds \Big) \\
& \hspace{0.5cm}- F \Big(t, p(t), \int_{g_1}^{j_1}
...\int_{g_m}^{j_m}
K(t, s, p(s)) ds, \int_{g_1}^{h_1} ... \int_{g_m}^{h_m} H(t, s, p(s)) ds \Big) \Big|\\
&\leq (1 - \delta_n)| a_n(t) - p(t)| + \delta_n \Big (\alpha |a_n(t) - p(t)|\\
& \hspace{0.5cm}+ \beta \Big |\int_{g_1}^{j_1} ... \int_{g_m}^{j_m}
K(t, s, a_n(s))ds - \int_{g_1}^{j_1} ... \int_{g_m}^{j_m} K(t, s,
p(s))
ds \Big| \\
& \hspace{0.5cm} + \gamma \Big |\int_{g_1}^{h_1} ...
\int_{g_m}^{h_m} H(t, s, a_n(s)) ds - \int_{g_1}^{h_1} ...
\int_{g_m}^{h_m} H(t, s,
p(s)) ds \Big| \Big)\\
& \leq (1 - \delta_n)| a_n(t) - p(t)| + \delta_n \alpha |a_n(t) -
p(t)| + \delta_n \beta \int_{g_1}^{j_1} ... \int_{g_m}^{j_m} L_K
|a_n(s) - p(s)| ds \\
& \hspace{0.5cm}+ \delta_n \gamma \int_{g_1}^{h_1} ...
\int_{g_m}^{h_m} L_H |a_n(s) - p(s)|ds\\
&\leq  \Big \{ 1 - \delta_n \Big (1 - \Big [\alpha +( \beta L_K +
\gamma L_H)\prod\limits_{i=1}^{m}(h_i - g_i) \Big] \Big ) \Big \}
\|a_n - p\|. \hspace{3cm}(5.4)
\end{aligned}
\end{equation*}
\begin{equation*}
\begin{aligned}
\|a_n - p\| &= \|Ac_n - p\| \\
&= |A(c_n)(t) - A(p)(t)|\\
 &= \Big| F \Big(t, c_n(t), \int_{g_1}^{j_1} ... \int_{g_m}^{j_m} K(t, s,
c_n(s)) ds, \int_{g_1}^{h_1} ... \int_{g_m}^{h_m} H(t, s, c_n(s))ds \Big)\\
& \hspace{0.5cm}-  F \Big(t, p(t), \int_{g_1}^{j_1}
...\int_{g_m}^{j_m}
K(t, s, p(s)) ds, \int_{g_1}^{h_1} ... \int_{g_m}^{h_m} H(t, s,p(s))ds \Big) \Big|\\
&\leq \Big [ \alpha +(\beta L_K + \gamma L_H)
\prod\limits_{i=1}^{m}(h_i - g_i) \Big ] \|c_n - p\|.
\hspace{5cm}(5.5)
\end{aligned}
\end{equation*}
Combining $(5.3)$, $(5.4)$ and $(5.5)$, we get
\begin{equation*}
\begin{aligned}
\| c_{n+1} - p\| & \leq \Big [ \alpha +(\beta L_K + \gamma L_H)
\prod\limits_{i=1}^{m}(h_i - g_i) \Big] \\
& \hspace{0.5cm} \times \Big \{ 1 - \delta_n \Big (1 - \Big [\alpha
+( \beta L_K + \gamma L_H)\prod\limits_{i=1}^{m}(h_i - g_i) \Big] \Big) \Big \}  \\
& \hspace{0.5cm} \times \Big [ \alpha +(\beta L_K + \gamma L_H)
\prod\limits_{i=1}^{m}(h_i - g_i) \Big] \times \|c_n - p\|\\
& \leq \Big \{ 1 - \delta_n \Big (1 - \Big [\alpha +( \beta L_K +
\gamma L_H)\prod\limits_{i=1}^{m}(h_i - g_i) \Big] \Big ) \Big \}
\times \|c_n - p\|.
\end{aligned}
\end{equation*}
Thus, by induction, we get
\begin{equation*}
\begin{aligned}
\| c_{n+1} - p\| & \leq  \|c_0 - p\| \times\\
& \hspace{0.5cm} \prod\limits_{k=0}^{n} \Big \{ 1 - \delta_k \Big(1
- \Big [\alpha +( \beta L_K + \gamma L_H)\prod\limits_{i=1}^{m}(h_i
- g_i)\Big ] \Big ) \Big \}. \hspace{3cm}(5.6)
\end{aligned}
\end{equation*}
Since $\delta_k \in [0, 1]$ for all $k \in \mathbb{N},$ assumption
$(A_5)$ yields\\
$$1 - \delta_k \Big (1 - \Big [\alpha +( \beta L_K +
\gamma L_H)\prod\limits_{i=1}^{m}(h_i - g_i)\Big ] \Big)  < 1.$$
Using the fact that $e^x \geq 1 - x$ for all $x \in [0, 1],$ we have
$$\| c_{n+1} - p\| \leq  \|c_0 - p\|e^{-(1 - [\alpha +( \beta L_K +
\gamma L_H)\prod\limits_{i=1}^{m}(h_i - g_i)]) \sum_{k=0}^{n}
\delta_k}$$ which yields $\lim\limits_{n \to\infty}\| c_n - p\| =
0.$

\section{ACKNOWLEDGEMENTS}
The authors are grateful to University Grants Commission, India for
providing financial assistance in the form of the BSR Start-Up
Research Grant and Junior Research Fellowship.

\end{document}